\documentclass[secthm,seceqn,amsthm,ussrhead,12pt]{amsart}
\usepackage{amsmath,latexsym}
\usepackage[english]{babel}
\usepackage[psamsfonts]{amssymb}
\usepackage{times}
\usepackage{cite}

\usepackage[mathcal]{euscript}

\newtheorem{Th}{Theorem}
\newtheorem{Lem}[Th]{Lemma}
\newtheorem{proposition}[Th]{Proposition}
\newtheorem{Remark}[Th]{Remark}

\newtheorem{Def}[Th]{Definition}

\begin{document}
\sloppy


\begin{center}

\textbf{A characterization of nilpotent nonassociative algebras\\ by invertible Leibniz-derivations.}

\medskip

\medskip
\textbf{Ivan Kaygorodov, Yury Popov}

\medskip

\sloppy
\sloppy

\medskip


\medskip

\section*{Abstract}
\end{center}

\medskip

Moens proved that a finite-dimensional Lie algebra over a field of characteristic zero is nilpotent if and only if it has an invertible Leibniz-derivation.
In this article we prove the analogous results for
finite-dimensional
Malcev, Jordan, $(-1,1)$-, right alternative, Zinbiel and Malcev-admissible noncommutative Jordan  algebras over a field of characteristic zero.
Also, we describe all Leibniz-derivations of semisimple Jordan, right alternative and Malcev algebras.

\medskip

{\bf Keywords:} Leibniz derivation, Malcev algebra, Jordan algebra, $(-1,1)$-algebra, nilpotent algebra.

\medskip

\section{Introduction}

\medskip

The theory of Lie algebras having a nonsingular derivation has a rich history, and is still an active research area.
Such Lie algebras appear in many different situations, such as in studies of pro-$p$ groups of finite coclass by
Shalev and Zelmanov \cite{Shalev94,SZ92} and in the problems concerning the existence of left-invariant affine structures on Lie groups (see Burde's survey \cite{Burde06} for details).

In 1955, Jacobson \cite{Jac} proved that a finite-dimensional Lie algebra over a
field of characteristic zero admitting a non-singular (invertible)
derivation is nilpotent. The problem of whether the inverse of this
statement is correct remained open until work \cite{Dix}, where
an example of a nilpotent Lie algebra in which all derivations are
nilpotent (and hence, singular) was constructed.
For Lie algebras in prime characteristic the situation is more complicated. In that
case there exist non-nilpotent Lie algebras, even simple ones, which admit nonsingular derivations \cite{BKK95}.
The main examples of nonsingular derivations are periodic derivations.
Kostrikin and Kuznetsov \cite{KK} noted that a Lie algebra admitting a nondgenerate derivation admits a periodic derivation, that is, a derivation $d$ such that $d^N = id$ for some $N$, and proved that a Lie algebra admitting a derivation of period $N$ is abelian provided that $N \not \equiv 0 \text{ (mod 6)}$.
Burde and Moens proved that a finite-dimensional complex Lie algebra L admits a periodic derivation if and only if L admits a nonsingular derivation whose inverse is again a derivation if and only if L is hexagonally graded \cite{BurdeMoens12}.
In the case of positive characteristic $p$, Shalev proved that if a Lie algebra $L$
admits a nonsingular derivation of order $n=p^sm,$ where $(m,p)=1$ and $m<p^2-1,$
then $L$ is nilpotent \cite{Shalev99}.
The study of periodic  derivations was continued by Mattarei \cite{mat02,mat07,mat09}.

The study of derivations of Lie algebras lead to the appearance of the notion of their
natural generalization --- a pre-derivation of a Lie algebra, which is a derivation of a Lie triple system induced by that algebra.
In \cite{Baj} it was proved that Jacobson's result is
also true in the terms of pre-derivations.
Several examples of nilpotent Lie
algebras which pre-derivations are nilpotent were presented in
\cite{Baj,Burd}.

In paper \cite{Moens} a generalization of derivations and
pre-derivations of Lie algebras is defined as a Leibniz-derivation of
order $k$.
Moens proved that a finite-dimensional Lie algebra over a field of characteristic zero is nilpotent if and only if
it admits an invertible Leibniz-derivation.
After that,
Fialowski, Khudoyberdiyev and Omirov \cite{Fialc12} showed
that with the definition of Leibniz-derivations from \cite{Moens}
the similar result for non-Lie Leibniz algebras is not true.
Namely, they gave an example of non-nilpotent Leibniz algebra which
admits an invertible Leibniz-derivation. In order to extend the
results of the paper \cite{Moens} for Leibniz algebras they introduced the
definition of Leibniz-derivation of a Leibniz algebra which agrees
with the definition of Leibniz-derivation of a Lie algebra and proved
that a finite-dimensional Leibniz algebra is nilpotent if and only if it admits an
invertible Leibniz-derivation.
It should be noted that there exist non-nilpotent Filippov ($n$-Lie) algebras with invertible derivations (see \cite{Williams}).
In paper \cite{KP} authors showed that the same result holds for alternative algebras (particularly, for associative algebras).
Also, in this article an example of nilpotent alternative (non-associative) algebra over a field of positive characteristic possessing only singular derivations was provided.

It is well known that the radicals of finite-dimensional algebras belonging to the classical varieties (such as varities of Jordan algebras, Lie algebras, alternative algebras, and many others) are invariant under their derivations \cite{Slinko1972}. Therefore, it is natural to state another interesting problem: is the radical of an algebra invariant under its Leibniz-derivations? Moens proved that the solvable radical of a Lie algebra is invariant under all its Leibniz-derivations, Fialowski, Khudoyberdiev and Omirov showed the invariance of solvable and nilpotent radicals of a Leibniz algebra, and the authors proved analogous results for alternative algebras.
Another interesting task is to describe the Leibniz-derivations of algebras belonging to certain "nice" classes, such as semisimple and perfect algebras. Moens described all Leibniz-derivations of semisimple Lie algebras \cite{Moens}
and Zhou described all pre-derivations of perfect centerless Lie algebras \cite{Zhou} of characteristic $\neq 2$.

Nowadays, studies of nonassociative algebras with derivations draw a lot of attention,
in particular, Jordan and Malcev algebras.
For example,
A. Popov determined the structure of differentiably simple Jordan algebras \cite{Popov2},
Kaygorodov and A. Shestakov studied generalized derivations of Jordan algebras  \cite{kay,kay_lie2,shestakov},
Filippov and Kaygorodov described generalized derivations of Malcev algebras \cite{Filll,Kay_Mal}.

\medskip

The main purpose of this article is to prove the analogues of
Moens' theorem for Jordan, $(-1,1)$- and Malcev algebras, that is, to prove that an algebra belonging to one of this classes is nilpotent if and only if it admits an invertible Leibniz-derivation. Also, we prove that any Leibniz-derivation of semisimple Jordan, $(-1,1)$- or Malcev algebra is a derivation of this algebra and show that the radicals of algebras of these classes are invariant with respect to Leibniz-derivations. We also discuss the relations between nilpotency and the property of having invertible Leibniz-derivations for many classes of nonassociative algebras, such as right alternative, quasiassociative, quasialternative, noncommutative Jordan Malcev-admissible, $n$-Lie (Filippov) and $n$-ary Malcev algebras. We give some sufficient conditions under which these algebras are nilpotent.

This paper is organized as follows. In Section 2, we discuss and present the notion of Leibniz-derivation, state some preliminary lemmas and discuss some of the previous results on the topic. In Section 3 we prove that the right nilpotent radical of an algebra is invariant under its left Leibniz-derivations provided that some additional conditions hold. In Section 4 we prove our main results for Malcev algebras. Namely, we prove that the solvable radical of a Malcev algebra is invariant under its left Leibniz-derivations, prove that for a semisimple Malcev algebra the sets of its left Leibniz-derivations and derivations coincide, and finally prove that a finite-dimensional Malcev algebra over a field of characteristic zero is nilpotent if and only if it is has an invertible Leibniz-derivation. In Section 5 we prove the Moens' theorem for Jordan and $(-1,1)$-algebras. Also we find sufficient conditions for nilpotency of noncommutative Jordan Malcev-admissible and quasialternative algebras.
In Section 6 we discuss some remarks and open questions.

Unless otherwise stated, 
in what follows all spaces of algebras are assumed finite-dimensional over the field of characteristic zero.

\medskip

\section{Preliminaries}

\medskip

We are using the standard notations:

$$(x,y,z):=(xy)z-x(yz)\mbox{ --- the associator of elements }x, y, z,$$
$$\left[x,y\right]:=xy-yx\mbox{ --- the commutator of elements }x, y,$$

For an algebra $A$ we define the following chains of subsets:
\begin{eqnarray*} A^1 = A, A^n = A^{n-1}A + A^{n-2}A^2 + \ldots +AA^{n-1}; \end{eqnarray*}
\begin{eqnarray*} A^{\langle 1 \rangle} = A, A^{\langle n \rangle} = A^{\langle n-1 \rangle}A; \end{eqnarray*}
\begin{eqnarray*} A^{(1)} = A^2, A^{(n)} = (A^{(n-1)})^2. \end{eqnarray*}

The subset $A^n$ is called the \textit{n-th power of A}, the subset $A^{\langle n \rangle}$ is called the \textit{n-th right power of A}, and the subset $A^{(n)}$ is called the \textit{n-th solvable power of A}. An algebra $A$ is called \textit{nilpotent} (respectively, \textit{right nilpotent}, \textit{solvable}), if $A^n = 0$ (respectively, $A^{\langle n \rangle} = 0,$ $A^{(n)} = 0$) for some $n \in \mathbb{N}$ . The minimal such $n$ is called the \textit{index of nilpotency of} $A$ (\textit{index of right nilpotency of} $A$, \textit{index of solvability of} $A$).

It is known that for commutative and anticommutative algebras the notions of nilpotency and right nilpotency coincide:
\begin{proposition}\cite[p.95]{kolca}
\label{right_nilp}
Let $A$ be a commutative (anticommutative) algebra. Then $A^{2^n} \subseteq A^{\langle n \rangle}.$ In other words, if $A$ is right nilpotent of index $n$, then it is nilpotent of index not greater than $2^n.$
\end{proposition}

For an arbitrary element $x \in A$ consider the operators $L_x, R_x: A \rightarrow A$ of respectively left and right multiplication by $x$  defined by $L_x(a) = xa, R_x(a) = ax.$
The algebra $M(A) \subseteq End(A)$ generated by $L_x, R_x, x \in A$ is called the\emph{ multiplication algebra of $A.$} It is widely known that $A$ is nilpotent if and only if $M(A)$ is nilpotent.
We call the \textit{Lie transformation algebra} $T(A)$ the smallest Lie subalgebra of $End(A)^{(-)}$ containing the set $L(A) + R(A)$ spanned by the operators $L_x, R_x, x \in A.$
A derivation $d \in Der(A)$ is called \textit{inner}, if it lies in $T(A)$, and \textit{outer} otherwise. It is a well known fact that the set $InDer(A)$ of inner derivations of $A$ is an ideal of the algebra $Der(A).$

Now we consider the definitions of a Leibniz-derivation given in various sources. The first definition of a Leibniz-derivation is due to Moens \cite{Moens}:

\begin{Def}\label{defmoens}
\emph{A Leibniz-derivation (by Moens) of order $n$ of a Lie algebra $L$ is an endomorphism $d$
of that algebra satisfying the relation}
\begin{equation}
\label{moens_lder}
d(x_1(x_2 \dots (x_{n-1}x_n)\dots))=\sum\limits_{i=1}^n x_1(x_2 \dots d(x_i) \dots (x_{n-1}x_n)\dots).
\end{equation}
\end{Def}

In the paper \cite{Fialc12} it was shown that this definition is not substantial for the case of general right Leibniz algebras. In order to carry over the results of \cite{Moens} authors had to redefine the notion of a Leibniz-derivation according to their needs:

\begin{Def}\label{deffok}
\emph{A Leibniz-derivation (by Fialowski-Khudoyberdiev-Omirov) of order $n$ of a right Leibniz algebra $L$ is an endomorphism $d$
of that algebra satisfying the relation}
\begin{equation}
\label{FKO_lder}
d(((\dots(x_1x_2) \dots )x_{n-1})x_n)=\sum\limits_{i=1}^n ((\dots(x_1x_2) \dots d(x_i) \dots )x_{n-1})x_n.
\end{equation}
\end{Def}

It is easy to see that for  Lie algebras definitions by Moens and by Fialowski-Khudoyberdiev-Omirov agree.

Throughout the paper we work with many different classes of nonassociative algebras, hence we will need a more general definition of Leibniz-derivation which generalizes both (\ref{moens_lder}) and (\ref{FKO_lder}).
Firstly, let $f$ be an arrangement of brackets on a product of length $n$. The product of elements $x_1, \dots, x_n$ with the arrangement $f$ will be denoted by $[x_1, \dots, x_n]_f.$

\begin{Def} Let $A$ be an algebra, $n$ be a natural number $ \geq 2$, and $f$ be an arrangement of brackets on a product of length $n$. A linear mapping $d$ on $A$ is called an $f$-Leibniz derivation of $A$, if for any $a_1, \dots a_n \in A$ we have
\begin{equation}
\label{f_lder}
d([a_1,\ldots,a_n]_f) = \sum\limits_{i=1}^n[a_1,a_2,\ldots, d(a_i), \ldots a_n]_f.
\end{equation}
\end{Def}

Particularly, if $f = l(n)$ ($f = r(n)$) is the left (right) arrangement of brackets of length $n$, that is,
$$[x_1, \ldots, x_n]_{l(n)} = ((\ldots(x_1x_2) \ldots )x_{n-1})x_n, \
 \ [x_1, \ldots, x_n]_{r(n)} = (x_1(x_2 \ldots (x_{n-1}x_n)\ldots)) \ ,$$
then $l(n)$-Leibniz derivation ($r(n)$-Leibniz derivation) of $A$ will be called a \emph{left (right) Leibniz-derivation of $A$ of order $n$.}
If $d$ is an $f$-derivation for any arrangement $f$ of length $n$, then $d$ will be called a \emph{Leibniz-derivation of $A$ of order $n.$} One can see that in our terms a Leibniz-derivation by Moens is a right Leibniz-derivation of a Lie algebra $L$, and a Leibniz-derivation by Fialowski-Khudoyberdiev-Omirov is a left Leibniz-derivation of a Leibniz algebra $L.$
It is easy to see that for

(1) (anti)commutative algebras the notion of a left Leibniz-derivation coincides with the notion of a right Leibniz derivation;

(2) left Leibniz algebras every right Leibniz-derivation is a left Leibniz-derivation;

(3) left Zinbiel algebras every left Leibniz-derivation is a right Leibniz-derivation;

(4) associative 
and Lie algebras the notion of a Leibniz-derivation coincides with the notion of an $f$-Leibniz derivation for any $f$.

\medskip

This definition can be naturally interpreted using the language of $n$-ary algebras. For an algebra $A$ we can define on its underlying vector space
the structure of $n$-ary algebra $A_f$ with the multiplication
$$[a_1, a_2,\ldots, a_n]_n=[a_1,a_2, \dots, a_n]_f.$$
Similarly to the case of binary algebras, we say that a linear map $d$ is a derivation of $n$-ary algebra $B$ if $d$ satisfies the relation
$$d([x_1,\dots,x_n]_n) = \sum \limits_{i=1}^n [x_1, \dots, d(x_i), \dots, x_n]_n. $$
Then $d$ is an $f$-Leibniz-derivation of $A$ if and only if it is a derivation of $A_f$.

The notion of an ideal can be naturally generalized for $n$-ary multiplication: a subspace $I$ of an $n$-ary algebra $B$ is called an ideal if
$$[I,B,\dots ,B]_n \subseteq I, [B,I,\dots, B]_n \subseteq I, \dots, [B, \dots ,B ,I]_n \subseteq I.$$
An $n$-ary algebra $B$ is called \textit{n-solvable} if its derived sequence
$$ \ B^{(1)_n} = [B,\dots,B]_n, \ B^{(t+1)_n} = [B^{(t)_n},\dots,B^{(t)_n}]_n, \text{for all t} \geq 1$$
descends to zero.
For an $n$-ary algebra $C$ by $Rad(C)$ we denote the largest $n$-solvable ideal of $C$.
This radical is well defined since the sum of arbitrary two solvable ideals of $C$ is again a solvable ideal of $C$.

In this article we will mostly be working with left Leibniz-derivations, therefore, unless otherwise stated, for an algebra $A$ and a natural number $n \geq 2$ by $A_n$ we mean the $n$-ary algebra on the vector space $A$ with the product
$$[a_1, \dots, a_n]_n = [a_1, \dots, a_n]_{l(n)}.$$
\medskip

Let $A$ be an algebra and $f$ be an arrangement of brackets. By $LDer_f(A)$ we denote the space of all $f$-Leibniz-derivations of $A$. Also let $LDer_n(A)$ denote the space of all Leibniz-derivations of order $n$ of algebra $A$, $LDer(A) = \bigcup_{n \geq 2} LDer_n(A)$ be the set of all Leibniz-derivations of $A,$ and $LDer_l(A) = \bigcup_{n \geq 2} LDer_{l(n)}(A)$ denote the set of all left Leibniz-derivations of $A.$
Note that a (left) Leibniz-derivation of order 2 is a usual derivation of $A$. Moreover, trivial induction can show that any derivation of $A$ is a Leibniz-derivation of $A$ of any order $n \geq 2.$ Thus, the order of Leibniz-derivation is not unique.
Moreover, Moens proved the following Lemma (actually, he proved it only for the right Leibniz-derivations, but one can easily that this lemma is true for the left Leibniz-derivations and Leibniz-derivations):
\begin{Lem}\cite{Moens}
\label{ord_lem}
If $s, t \in \mathbb{N}$ and $s | t$ then $LDer_{l(s+1)}(A) \subseteq LDer_{l(t+1)}(A), LDer_{s+1}(A) \subseteq LDer_{t+1}(A).$
\end{Lem}

We will need the following proposition by Moens:

\begin{proposition}\cite{Moens}
\label{subalg_prop}
For $n \geq 2$ the subsets $LDer(A)$, $LDer_n(A)$ and $LDer_{l(n)}(A)$ are in fact subalgebras of the Lie algebra of endomorphisms of $A$ and we have the chains
$$ InDer(A) \subseteq Der(A) \subseteq LDer_n(A) \subseteq LDer(A) \subseteq End(A)$$
$$ InDer(A) \subseteq Der(A) \subseteq LDer_{l(n)}(A) \subseteq LDer_l(A) \subseteq End(A)$$
\end{proposition}

These chains collapse for algebras with special properties. For example, if an algebra $A$ satisfies the identity $[x_1,\dots,x_n]_f = 0$ for an arrangement of brackets $f$, then it is easy to see that any endomorphism of $A$ is in fact a $f$-Leibniz derivation of $A$, that is, $LDer_f(A) = LDer(A) = End(A).$ Actually, the converse of this statement is also true. Suppose that an algebra $A$ over the field of characteristic 0 satisfies the property $LDer_f(A) = End(A).$ Therefore, an identity map of $A$ is also an $f$-Leibniz-derivation, hence for any $a_1,\ldots,a_n \in A$ we have

\begin{equation*} [a_1,\dots,a_n]_f = \sum_{i = 1}^n [a_1,\dots,a_n]_f, \end{equation*}

and

\begin{equation*} (n-1)[a_1,\dots,a_n]_f = 0 \end{equation*}

and $A$ satisfies the identity $[x_1,\dots,x_n]_f = 0.$ 

Particularly, if we take $f = l(n)$ or assume that the above identity holds for all arrangements of brackets of length $n,$ we have
\begin{proposition}
\label{right_nilp_der}
An algebra $A$ is right nilpotent if and only if $LDer_{l(n)}(A) = LDer_l(A) = End(A)$ for some $n \geq 2.$ Moreover, the minimal $n$ with such property is equal to the index of right nilpotency of $A.$
\end{proposition}

\medskip

\begin{proposition}\label{nilp_der}
An algebra $A$ is nilpotent if and only if $LDer_n(A) = LDer(A) = End(A)$ for some $n \geq 2.$ Moreover, the minimal $n$ with such property is equal to the index of nilpotency of $A.$
\end{proposition}

It is also easy to describe the Leibniz-derivations of unital algebras:
\begin{proposition}
If $A$ is an algebra with unity $1$ over the field of characteristic zero then for any arrangement of brackets $f$ we have $LDer_f(A) = Der(A).$
\end{proposition}
Indeed, let $f$ be an arrangement of brackets on the product of $n$ elements. Then for any $d \in LDer_f(A)$ we have $d(1) = n d(1),$ so $d(1) = 0.$ Now, for $x, y \in A$

\begin{equation*} d(xy) = d([x,y,1,\dots,1]_f) = d(x)y + xd(y), \end{equation*}
therefore $d \in Der(A),$ and the first chain takes the following form:

\begin{equation*} InDer(A) \subseteq Der(A) = LDer_n(A) = LDer(A) \subseteq End(A).\end{equation*}

In particular, all $f$-Leibniz-derivations of semisimple  associative, Jordan, alternative, $(-1,1)$,
or right alternative algebra $A$ over the field of characteristic zero are derivations of $A$, even inner derivations.

\medskip

\begin{Lem}
\label{invert_der_exist}
Every nilpotent algebra of nilpotency index $n$ admits an invertible Leibniz-derivation of order $\lceil\frac{n}{2}\rceil.$
\end{Lem}
\textbf{Proof.} Let $A$ be an algebra with nilpotency index $n.$ Set $q = \lceil\frac{n}{2}\rceil$ and consider $A$ as the direct sum of vector spaces:
$A = A^q \oplus W.$ Now, we can define the linear mapping $\phi$ in the following way:
$$\phi(x)=\left\{\begin{array}{ll}
x,  \mbox{ if } x \in W,\\
nx, \mbox{ if } x \in A^n.
\end{array} \right.$$

It is easy to see that $\phi$ is a Leibniz-derivation of $A$ of order $\lceil\frac{n}{2}\rceil$.
The Lemma is proved.

\medskip

We will also need the following generalization of Leibniz rule for the left Leibniz--derivations of order $n:$
\begin{Lem}
\label{Leib_rule}
Let $A$ be an algebra and $d$ be a left Leibniz-derivation of order $n$ of $A$. Then the following relation holds for any $k \in \mathbb{N}:$
$$d^k([x_1,x_2,\ldots,x_n]_n)=\sum\limits_{i_1 + \ldots + i_n = k}\frac{k!}{i_1!\dots i_n!}[d^{i_1}(x_1), \ldots,d^{i_n}(x_n)]_n.$$
\end{Lem}

We remark that throughout this article without loss of generality we may assume that the base field $F$ is algebraically closed. Indeed, if it is not the case, we can extend $d$ to the algebra $\overline{A} = A \otimes_F \overline{F}$ over the algebraic closure $\overline{F}$ of $F$ by the rule $\overline{d}(x \otimes \alpha) = d(x) \otimes \alpha.$ It is easy to see that $\overline{d}$ is a Leibniz-derivation of $\overline{A}$ of the same order, and the map $\overline{d}$ is invertible if and only if $d$ is. Also, we can identify $A$ with the $F$-subalgebra $A \otimes 1$ of $\overline{A},$ hence $A$ is nilpotent if $\overline{A}$ is nilpotent.


\section{The invariance of the right nilpotent radical}

First, we prove that under some conditions all left Leibniz-derivations of an algebra preserve its right nilpotent radical.
\begin{Th}
\label{rad_inv}
Let $A$ be an algebra over a field of characterstic 0 and $Rad(A)$ be the right nilpotent radical of $A$. Suppose that\\
$1)$ $A = Rad(A) + S$, where $S \cong A/Rad(A)$ is the semisimple part of $A$;\\
$2)$ $S = S^2$;\\
$3)$ $S = A / Rad(A)$ contains no nonzero right nilpotent right ideals.\\
Then $Rad(A)$ is invariant under all left Leibniz-derivations of $A.$
\end{Th}
\textbf{Proof.}
Let $d$ be a left Leibniz-derivation of order $n$ of $A.$ It suffices to show that $(Rad(A) + d(Rad(A)))/Rad(A)$ is a right nilpotent right ideal of $A/Rad(A) = S$. Using the decomposition $A = Rad(A) + S$ and the fact that $S = S^2$, we have
$$A(Rad(A) + d(Rad(A))) = $$
$$ = (Rad(A) + S)(Rad(A) + d(Rad(A))) \subseteq Rad(A) + S d(Rad(A)) = $$
$$ = Rad(A) + S^{\langle n-1 \rangle} d(Rad(A)) \subseteq $$
$$\subseteq Rad(A) + d([S, S, \dots, S, Rad(A)]_n) + \sum [S, S, \dots, d(S), \dots, S, Rad(A)]_n \subseteq$$
$$\subseteq Rad(A) + d(Rad(A)),$$
so $Rad(A) + d(Rad(A))$ is a right ideal of $A.$
We now have to prove the right nilpotency of $Rad(A) + d(Rad(A))$, that is, we have to find a natural number $n$ such that $(Rad(A) +d(Rad(A)))^{\langle n \rangle} = 0.$ Let $n_0$ denote the index of right nilpotency of $Rad(A).$ Lemma \ref{ord_lem} implies that we may assume that $n \geq n_0.$ Using Lemma \ref{Leib_rule}, for any $a_1, a_2, \dots, a_n \in Rad(A)$ we have
$$[d(a_1),d(a_2),\dots, d(a_n)]_n =$$
$$\frac{1}{n!}\left( d^n([a_1,a_2,\dots,a_n]_n) -
\sum \limits_{\substack{i_1 + l\dots + i_n = n \\i_j = 0 \text{ for some } j}}\frac{n!}{i_1!\dots i_n!}[d^{i_1}(a_1),\dots,d^{i_n}(a_n)]_n \right) = $$
$$ = -\sum \limits_{\substack{i_1 + \ldots + i_n = n \\ i_j = 0 \text{ for some } j}}\frac{1}{i_1!\dots i_n!}[d^{i_1}(a_1),\dots,d^{i_n}(a_n)]_n \in Rad(J).$$

Therefore, $(Rad(A) + d(Rad(A)))^{\langle n \rangle} \subseteq Rad(A) + d(Rad(A))^{\langle n \rangle} \subseteq Rad(A).$ Consequently, $(Rad(A) + d(Rad(A)))/Rad(A)$ is a right nilpotent right ideal of $A/Rad(A) = S.$ From the condition $3)$ it follows that $d(Rad(A)) \subseteq Rad(A).$ Theorem is now proved.
\medskip

\section{Moens' theorem for Malcev algebras}

\medskip

An algebra $M$ with multiplication $[ \cdot , \cdot ]$ is called \textit{Malcev} if it satisfies the following identities:
\begin{eqnarray*} [x,x] = 0, J(x, y, [x,z]) = [J(x, y, z)x], \end{eqnarray*}
where $J(x, y, z) = [[x, y], z] + [[z, x], y] + [[y, z], x]$ is the Jacobian of elements $x, y, z$.
The class of Malcev algebras is a well-studied extension of the class of Lie algebras.
They first arose in \cite{Malcev} and have relations with alternative algebras and Moufang loops.
They were studied thoroughly in many articles, such as
\cite{Sagle,kornevshestakov} and many others.



In \cite{Kuzmin} it was shown that in an arbitrary Malcev algebra the following operator identity holds:
\begin{equation*} 2R_{[[x,y],z]} = [[R_x,R_y],R_z] + [R_y,R_{[z,x]}] + [R_x,R_{[y,z]}]. \end{equation*}

Using this relation, it is easy to compute the Lie multiplication algebra for any Malcev algebra $M$:
$$T(M) = R(M) + [R(M),R(M)].$$

In his work \cite{Sagle} Sagle has shown that an arbitrary inner derivation of a Malcev algebra $M$ is of form
\begin{equation*} \sum_i D(x_i,y_i) + R(n), \end{equation*}
where $x_i, y_i \in M$, $D(x,y) = [R_x,R_y] + R_{[x,y]}$, and $n$ is an element of $N(M)$ -- the \textit{Lie center} of $M$, which consists of elements $z$ such that $J(z,M,M) = J(M,z,M) = J(M,M,z) = 0.$

Using that fact, we can see that any element of $T(M)$, where $M$ is a Malcev algebra, is of the form $D + R_x$, where $D$ is an inner derivation of $M$, and $x$ is an element of $M.$ Also, if $N(M) = 0$, such representation is unique.

Kuzmin has shown that, like Lie algebras, Malcev algebras possess the \textit{Killing form} $\chi(x,y) = Tr(R_xR_y)$ which is symmetric and associative in a sense that $\chi([x,y],z) = \chi(x,[y,z])$ for all $x, y, z.$

The classical Engel and Lie theorems have the following analogues for Malcev algebras:
\begin{Th}\cite{Zhevl}
A Malcev algebra $M$ is nilpotent if and only if $R_x$ is nilpotent for any $x \in M.$
\end{Th}

\begin{Th}\cite{Kuzmin}
Let $\rho$ be a representation of a solvable Malcev algebra $M$ over the field $F$ such that for any $x \in M$ all eigenvalues of $\rho(x)$ lie in $F.$ Then all matrices $\rho(x)$ can be made upper-triangular.
\end{Th}

Is is a well-known fact that the  Malcev algebra also admits a Levi decomposition, that is, it can be regarded as the direct sum of its radical $Rad(M)$ and semisimple part $S$.

The classification of simple Malcev algebras over an algebraically closed field is widely known: aside from simple Lie algebras, there is only one non-Lie 7-dimensional Malcev simple algebra $M_7$, which can be obtained from the split Cayley--Dickson algebra by introducing a new operation $x*y = \frac{1}{2}[x,y]$, where $[,]$ is a usual commutator in Cayley--Dickson algebra, and taking a quotient algebra by the base field (for the full details see, for example, \cite{Kuzmin}).

We remind that for a Malcev algebra $M$ and $n \geq 2$ we denote by $M_n$  an $n$-ary algebra on the same vector space with multiplication $[x_1,x_2,\dots,x_n]_n = [[\dots[x_1,x_2],x_3,\dots],x_n].$

Our main goal for this paragraph is to prove that a Malcev algebra $M$ admitting a nonsingular left Leibniz-derivation of order $n$ is nilpotent.

We begin by proving that the solvable radical of $M$ and $n$-solvable radical of $M_n$ coincide:
\begin{Lem}
\label{radradn}
Let $M$ be a Malcev algebra and $n \geq 2.$ Then $Rad(M) = Rad(M_n).$
\end{Lem}

\medskip
\textbf{Proof.}
Easy induction can show that for any subalgebra $B \subseteq M$ its $n$-ary solvable power is a subset of its usual solvable power, i.e., $B_n^{(k)_n} \subseteq B^{(k)}$ for any $k \geq 1.$ Particulary, any solvable ideal of $M$ is a $n$-solvable ideal of $M_n,$ therefore $Rad(M) \subseteq Rad(M_n).$ Let $M = S \oplus Rad(M)$ be a Levi decomposition of $M$ and $\pi: M \rightarrow S$ be the canonical projection map. Since $\pi$ is a morphism of $M$, we have
$$ [S,\pi(Rad(M_n))] = [S^{\langle n-1 \rangle}, \pi(Rad(M_n))] = [S,\dots,S,\pi(Rad(M_n))]_n =$$
 $$= [\pi(M),\dots,\pi(M),\pi(Rad(M_n))]_n = \pi([M,\dots,M,Rad(M_n)]_n) \subseteq \pi(Rad(M_n)), $$
therefore $\pi(Rad(M_n))$ is an ideal of $S$, therefore, is perfect: $\pi(Rad(M_n))^2 = \pi(Rad(M_n)).$ But since $\pi$ is a homomorphism, $\pi(Rad(M_n))$ must be an $n$-solvable $n$-ary algebra. Only zero $n$-ary algebra can be both $n$-solvable and perfect, therefore $\pi(Rad(M_n)) = 0$ and $Rad(M_n) \subseteq Rad(M).$ The Lemma is now proved.

\medskip

\begin{Lem}
Let $M$ be a Malcev algebra, $I$ be an ideal of $M_n$, and $d \in LDer_n(M)$ be a left Leibniz-derivation of order $n$ of $M.$ Then
$$d(I)^{(k)_n} \subseteq I + d^{n^k}(I^{(k)_n})$$
for any $k \in \mathbb{N}.$
\end{Lem}

\medskip

\textbf{Proof.} We prove this statement by induction on $k.$ For $k = 1,$ by Lemma \ref{Leib_rule} we have
$$d(I)^{(1)_n} = [d(I),\dots,d(I)]_n = d^n([I,\dots,I]_n) +$$
 $$ + \sum \limits_{\substack{i_1 + i_2 + \ldots + i_n = n \\i_j = 0 \text{ for some } j}} [d^{i_1}(I),\dots,d^{i_{j-1}}(I),I,d^{i_{j+1}}(I),\dots,d^{i_n}(I)]_n \subseteq $$
  $$\subseteq I + d^n(I^{(1)_n}).$$
  Assume now that the statement of the Lemma holds for $k \in \mathbb{N}.$ Using again Lemma \ref{Leib_rule},
  we verify the Lemma for $k+1:$
$$d(I)^{(k+1)_n} = [d(I)^{(k)_n},\dots,d(I)^{(k)_n}]_n \subseteq [I + d^{n^k}(I^{(k)_n}),\dots,I + d^{n^k}(I^{(k)_n})]_n \subseteq$$
$$ \subseteq I + d^{n^{k+1}}([I^{(k)_n},\dots,I^{(k)_n}]_n)$$
  $$\subseteq I + d^{n^{k+1}}(I).$$
  The Lemma is now proved.

\medskip

Now we prove that the solvable radical is invariant under all left Leibniz-derivations.

\medskip

\begin{Th}
\label{rad_invar}
Let $M$ be a Malcev algebra over a field of characteristic zero. 
Then $Rad(M)$ is invariant under all left Leibniz-derivations of $M.$
\end{Th}
\textbf{Proof.} From Proposition \ref{radradn},
it follows that we have to prove the invariance of $Rad(M_n).$ Let $k \in \mathbb{N}$ be such that $Rad(M_n)^{(k)_n} = 0.$ Previous Lemma implies that
$$d(Rad(M_n))^{(k)_n} \subseteq Rad(M_n) + d^{n^k}(Rad(M_n)^{(k)_n}) = Rad(M_n).$$
Therefore
$$(Rad(M_n) + d(Rad(M_n)))^{(2k)_n} = ((Rad(M_n) + d(Rad(M_n)))^{(k)_n})^{(k)_n} \subseteq$$
$$\subseteq (Rad(M_n) + d(Rad(M_n))^{(k)_n})^{(k)_n} = Rad(M_n)^{(k)_n} = 0.$$
It is also easy to verify that $Rad(M_n) + d(Rad(M_n))$ is an ideal of $M_n:$
$$[M, \dots,M,Rad(M_n) + d(Rad(M_n)),M,\dots,M]_n \subseteq$$
$$ \subseteq Rad(M_n) +  d([M,\dots,M,Rad(M_n),M,\dots,M]_n) +$$
$$+ \sum[M,\dots,d(M),M,\dots,M,Rad(M_n),M,\dots,M]_n +$$
$$+ \sum[M,\dots,M,Rad(M_n),M,\dots,M,d(M),M,\dots,M]_n \subseteq Rad(M_n) + d(Rad(M_n)).$$

Therefore $Rad(M_n) + d(Rad(M_n))$ is an $n$-solvable ideal of $M_n.$ Since $Rad(M_n)$ is the $n$-solvable radical of $M_n,$ it follows that $d(Rad(M_n)) \subseteq Rad(M_n).$ The theorem is now proved.

\medskip

We now show that semisimple Malcev algebras do not admit invertible left Leibniz-derivations.
In fact, we prove even stronger result:
\begin{Th}
\label{Leib_der_ss}
Let $M$ be a semisimple Malcev algebra over a field of characteristic 0. Then $Der(M) = LDer_l(M).$
\end{Th}

Since a semisimple  Malcev algebra $M$ is a direct
sum of simple Malcev algebras $M_i, i =1,\dots,n$ 
where $M_i = M_i^k$ for any $i = 1,\dots,n,~ k \in \mathbb{N},$
any left Leibniz-derivation of $M$ leaves each $M_i$ invariant,
so in fact we only have to prove the Theorem for simple Malcev algebras.
The analagous result for the semisimple Lie algebras was proved by Moens in \cite{Moens}.
Therefore we may assume that $M = M_7,$
though during the proof we provide notes showing how to adapt it for the case of any semisimple Malcev algebra.

We begin with the following Lemma, which generalizes the well known rule for commutating derivations and left multiplications:

\begin{Lem}
Let $A$ be an arbitrary algebra, $d$ be a left Leibniz-derivation of order $n$ of $A,$ $x_1,\dots, x_{n-1} \in A.$ Then
\begin{equation*} [d,L_{[x_1,\dots,x_{n-1}]_{n-1}}] = \sum_{i = 1}^{n-1}L_{[x_1,\dots,x_{i-1},d(x_i),x_{i+1},\dots,x_{n-1}]_{n-1}}. \end{equation*}
\end{Lem}
The proof of this Lemma is obvious.

\medskip

Since Malcev algebras are anticommutative, for arbitrary Malcev algebra we have the following relation:
\begin{equation}
\label{commutator} [d,R_{[x_1,\dots,x_{n-1}]}] = \sum_{i = 1}^{n-1}R_{[x_1,\dots,x_{i-1},d(x_i),x_{i+1},\dots,x_{n-1}]_{n-1}}.
\end{equation}

Now, for semisimple Malcev algebra $M$ we have $M = M^2,$ hence $M = M^{\langle n-1\rangle}$ and arbitrary element of $M$ can be represented as a sum of $(n-1)$-ary commutators:
\begin{equation*} x = \sum_i[x_{i_1},\dots,x_{i_{n-1}}]_{n-1}, x_{i_j} \in M. \end{equation*}

Applying (\ref{commutator}) to this relation, we see that

\begin{equation}
\label{ideal} [R(M),LDer_l(M)] \subseteq R(M).
\end{equation}

We are now able to prove
\begin{Lem}
Let $M$ be a semisimple Malcev algebra. The sum of spaces $T(M) + LDer_l(M) = L$ forms a Lie algebra under commutation of operators, and $T(M)$ is a semisimple ideal of $L.$
\end{Lem}

\medskip
\textbf{Proof.}
The fact that $LDer_l(M)$ is a Lie algebra follows from Proposition \ref{subalg_prop}.
Now, we know that $T(M) =  R(M) + [R(M),R(M)],$ so by (\ref{ideal}) and the Jacobi identity we have
$$[T(M),LDer_l(M)] = [R(M) + [R(M),R(M)], LDer_l(M)] \subseteq$$
$$\subseteq R(M) + [[R(M),LDer_l(M)],R(M)] + [R(M),[R(M),LDer_lR(M)]] \subseteq$$
$$\subseteq R(M) + [R(M),R(M)] = T(M),$$
so $T(M) \unlhd L.$ Semisimplicity of $T(M)$ follows from \cite{Kuzmin}. The Lemma is now proved.

\medskip

We want to study the structure of Lie algebra $L$ deeper:
\begin{Lem}
$T(M)$ is a direct summand of $L.$ Moreover, we have $L = T(M) \oplus T(M)^{\bot},$ where $T(M)^{\bot}$ is the orthogonal complement to $T(M)$ with respect to the Killing form on $L.$
\end{Lem}
\medskip
\textbf{Proof.} Since $T(M) \unlhd L$, the Killing form of $T(M)$ coincides with the restriction of $\chi_L$ to $T(M).$ Since $T(M)$ is a semisimple Lie subalgebra of $L,$ the restriction $\chi_{T(M)}$ of Killing form $\chi_L$ of $L$ to $T(M)$ is nondegenerate. Therefore, for any basis $x_1, \dots, x_k$ of $T(M)$ there exists a dual basis $y_1, \dots, y_k$ such that $\chi_{T(M)}(x_i,y_j) = \delta_{ij}.$ For arbitrary element $x \in L$ consider the element $x' = \sum_{i=1}^k\chi_L(x_i,x)y_i.$ For all $x_i, i=1,\dots,k$ we have $\chi_L(x_i,x-x') = 0,$ hence $x-x' \in T(M)^{\bot}$ and $L = T(M) + T(M)^{\bot}.$ Since $T(M)$ is a semisimple Lie algebra, $\chi_{T(M)}$ is nondegenerate, and this sum is direct. The Lemma is now proved.
\medskip

We are trying to prove that $L = T(M)$, so we need to show that $T(M)^{\bot} = 0.$
We proceed in a few steps:
\medskip
\begin{Lem}
\label{scalar_complement}
The orthogonal complement of $T(M)$ in $L$ consists only of scalar mappings: $T(M)^{\bot} \subseteq F \cdot id.$ Moreover, if $M$ is a simple Lie algebra, then $T(M)^{\bot} = 0.$
\end{Lem}
\textbf{Proof.} Since the orthogonal complement with respect to Killing form of an ideal is also an ideal, then $[T(M),T(M)^{\bot}] = 0.$ As before, we may assume that $M$ is a simple Malcev algebra. Let us suppose that $M = M_7.$ The structure of $T(M_7)$ was described in \cite{Kuzmin}: it was proved that
$T(M_7) \cong B_3$, the algebra of skew-symmetric transformations of $M_7$ with respect to its Killing form $\chi_{M_7}$. Since the base field $F$ is algebraically closed, there exists a basis of $M_7$ in which $\chi$ is represented by the identity matrix. Hence, in this basis the elements of $T(M_7)$ are represented by usual skew-symmetric matrices. It is easy to show that skew-symmetric matrices in fact generate the whole associative matrix algebra $End(M_7):$ using the standard notations $e_{ij}$ of matrix units, we see that for $i \neq j$ we have
\begin{equation*} (e_{ij} - e_{ji})(e_{ik} - e_{ki}) = -e_{jk}, \end{equation*}
\begin{equation*} (e_{ij} - e_{ji})^2 = -e_{ii} - e_{jj}. \end{equation*}
Using first relation, we can generate all $e_{ij}, i \neq j.$ Diagonal elements are easily obtained with the aid of second relation:
\begin{equation*} e_{ii} = \frac{1}{2}((e_{ii} + e_{11}) - (e_{11} + e_{jj}) + (e_{ii} + e_{jj})). \end{equation*}
Therefore, $End(M_7) = T(M_7)^2.$ Since for any $x \in End(M_7)$  the mapping $[x,\cdot]$ is a derivation of associative algebra $End(M_7)$ and $[x,T(M_7)] = 0$ for any $x \in T(M_7)^{\bot},$ by Leibniz rule we have
$$[x,End(M_7)] = [x,T(M_7)^2] = 0,$$
so $x$ lies in the center of the full matrix algebra $End(M_7)$, therefore, it is a scalar mapping.

In case when $M$ is a simple Lie algebra, the Jacobi identity implies that $T(M) = R(M) = Der(M),$ therefore $L = LDer_l(M).$ To prove the result of the Lemma, one may, for example, apply Schur's Lemma. In fact, we may even prove that a left Leibniz-derivation $d \in R(M)^{\bot}$ is zero mapping:
if $[d,R_x] = 0$ for any $x \in L$ then in anticommutative algebra we have
$[d,L_x] = 0$ and
$$ d([x_1,\dots,x_n]_n) =  [x_1,\dots,x_{i-1},d(x_i),x_{i+1},\dots,x_n]_n \text{ for any } i = 1, \dots, n.$$

But this means that
\begin{equation*} d([x_1,\dots,x_n]_n) = \sum_{i=1}^n [x_1,\dots,x_{i-1},d(x_i),x_{i+1},\dots,x_n]_n = n\cdot d([x_1,\dots,x_n]_n).  \end{equation*}
Since $L = L^{\langle n \rangle}, d = 0,$ and the proof of our theorem for semisimple Lie algebras finishes here. The Lemma is now proved.
\medskip

To prove the next Lemma, we will need the notion of ternary derivation and generalized derivation.

\begin{Def} Let $A$ be an algebra. The triple $(D,F,G) \in End(A)^3$  is called a \textit{ternary derivation} of $A$ if the following relation holds in $A:$
$$D(xy) = F(x)y + xG(y), x, y \in A.$$
The endomorphism $D$ of $A$ is called a \textit{generalized derivation} if there exist $F, G \in End(A)$ such that $(D,F,G)$ is a ternary derivation of $A.$
\end{Def}
Generalized and ternary derivations of algebras and superalgebras were studied by authors in \cite{kp16},
other authors in \cite{LL}, and in many other articles.

\begin{Lem}
The orthogonal complement of $T(M_7)$ in $L(M_7)$ is zero.
\end{Lem}
\medskip
\textbf{Proof.} From the work \cite{Kuzmin} it follows that the Lie center of $M_7$ equals zero, so $T(M_7) = Der(M_7) \oplus R(M_7).$ Hence, Lemma \ref{ord_lem} implies that any element of $T(M_7)^{\bot}$ has the form $d - R_x,$ where $d$ is a left Leibniz-derivation of $M_7,$ and $x \in M_7.$ From the previous Lemma we have
$d - R_x = \alpha \cdot id, \alpha \in F.$ Therefore, for arbitrary $y \in M_7$ from the relation (\ref{ideal}) we have
$$D(y,x) = R_{[y,x]} + [R_y,R_x] = R_{[y,x]} + [R_y, R_x + \alpha \cdot id] =$$
$$ = R_{[y,x]} + [R_y,d] \in R(M_7) \cap Der(M_7) = 0.$$
Applying the zero operator $D(x,y)$ on arbitrary element $z \in M_7,$ we get
\begin{equation*} [[y,z],x] = - [[y,x],z] + [y,[z,x]]. \end{equation*}
This means that the triple $(R_x,-R_x,R_x)$ is a ternary derivation of $M_7.$ From Kaygorodov's description \cite{Kay_Mal} of ternary and generalized derivations of $M_7$ it follows that
$R_x = D + \beta \cdot id,$ where $D \in Der(M_7), \beta \in F.$ Equating traces of both sides we infer that $\beta = 0.$ That means that $R_x \in R(M_7) \cap Der(M_7) = 0,$ and since $M_7$ is centerless, $x = 0$, and $d = \alpha \cdot id.$ It is easy to see that $\alpha = 0.$ For example, we can repeat verbatim the last part of the proof of previous Lemma: since $\alpha \cdot id$ is a scalar operator, $d$ commutes with all right multiplications, and so on. The Lemma is now proved.
\medskip

We are now ready to prove our theorem.

\textbf{Proof of  Theorem \ref{Leib_der_ss}.} Let $M$ be a semisimple Malcev algebra over the field $F$ of characteristic 0.
As before, we may assume that $M = M_7.$ By the previous Lemma we have $L = T(M_7) + LDer_l(M_7) = T(M_7),$ hence $LDer_l(M_7) \subseteq T(M_7).$ Since $T(M_7) = Der(M_7) \oplus R(M_7)$ we only have to prove that operators $R_x, x \in M_7$ cannot be left Leibniz-derivations of any order. This can be done by repeating verbatim the proof of previous Lemma with $\alpha = 0$. Therefore, we have $LDer_l(M_7) \cap R(M_7) = 0,$ and $LDer_l(M_7) = Der(M_7).$ The Theorem is now proved.

\medskip

\begin{proposition}
\label{Malcev_solv}
If a Malcev algebra $M$ admits an invertible left Leibniz-derivation, then $M$ is solvable.
\end{proposition}
\textbf{Proof.} As usual, let $Rad(M)$ be the solvable radical of $M$. From Theorem \ref{rad_invar} we see that $d(Rad(M)) \subseteq Rad(M),$ so $d|_S$ is an invertible left Leibniz-derivation of semisimple Levi component $S$ of $M.$ Suppose that $S \neq 0.$ Theorem \ref{Leib_der_ss} implies that $d|_S$ is a derivation of $S.$ We again assume that $S$ is a simple Malcev algebra. If $S$ is a simple Lie algebra, Jacobson's theorem implies that $S = 0$, a contradiction. If $S = M_7,$ the proof of Lemma \ref{scalar_complement} implies that $Der(M_7)$ is a subalgebra of the algebra $B_3$ of skew-symmetric $7 \times 7$ matrices, and it is widely known that a skew-symmetric matrix of an odd order is degenerate, so $M_7$ cannot have an invertible derivation. We have obtained a contradiction which implies that $S = 0$, and $M = Rad(M)$ is solvable. The Proposition is now proved.

\medskip

The nilpotency of $M$ will be established in a similar way to \cite{Moens}.

\medskip

\begin{Lem}
Let $M$ be a Malcev algebra, $d$ be its left Leibniz-derivation of order $n$, and $M = \bigoplus M_\alpha$ be its decomposition into sum of eigenspaces with respect to $d: M_\alpha = \bigcup_{k \geq 1} \ker (d - \alpha \cdot id)^k.$ Then for  eigenvalues $\alpha_1,\dots,\alpha_n$ we have
$$[M_{\alpha_1},\dots,M_{\alpha_n}]_n \subseteq M_{\alpha_1 + \ldots + \alpha_n}.$$
\end{Lem}
\textbf{Proof.} Using induction on $k \in \mathbb{N},$ one can similarly to Lemma \ref{Leib_rule} prove that the following relation holds in $M:$
$$
(d-(\alpha_1 + \ldots + \alpha_n)\cdot id)^k([x_1,\dots,x_n]_n)=$$
$$ = \sum\limits_{i_1 + \ldots + i_n = k}\frac{k!}{i_1!i_2!\dots i_n!} [(d-\alpha_1 \cdot id)^{i_1}(x_1), \dots,(d-\alpha_n \cdot id)^{i_n}(x_n)]_n.$$
By the definition of $M_{\alpha_i}$ there exist $k_i \in \mathbb{N}$ such that $(d - \alpha_i \cdot id)^{k_i}(M_{\alpha_i}) = 0.$ If we take $k = k_1 + \ldots + k_n$ in the previous equation, for arbitrary $x_1 \in M_{\alpha_1},\dots x_n \in M_{\alpha_n}$ we have
$$(d-(\alpha_1 + \ldots + \alpha_n)\cdot id)^k([x_1,\ldots,x_n]_n)= 0.$$
The Lemma is now proved.

\medskip

\begin{proposition}
\label{Malcev_nilp}
If $M$ is a solvable Malcev algebra that admits a invertible left Leibniz-derivation $d$ of order $n$, then $M$ is nilpotent.
\end{proposition}
\textbf{Proof.} By the previous Lemma, for arbitrary pair of eigenvalues $\alpha, \beta$ of $d$ we have $(R_{M_{\alpha}})^{n-1}(M_{\beta}) \subseteq M_{\beta + (n-1)\alpha}.$ Since the set of eigenvalues of $d$ is finite, it is possible to choose a $k \in \mathbb{N}$ such that $(R_{M_{\alpha}})^{(n-1)k}(M_{\beta}) = 0$ for arbitrary eigenvalues $\alpha, \beta$ of $d.$ Therefore, for any $x \in M_{\alpha}, ~ R_x$ is a nilpotent operator. Since $M$ is a direct sum of its eigenspaces, we may choose a basis of $M$ consisting entirely of vectors $x$ such that $R_x$ is nilpotent. Since $M$ is a solvable Malcev algebra, Lie's theorem implies that all operators $R_x, x \in M$ are nilpotent, and from Engel's theorem it follows that $M$ is nilpotent. The Proposition is now proved.

\medskip

Now we can prove the main theorem of this paragraph:

\begin{Th}
\label{Moens_th_Malc}
A Malcev algebra $M$ over a field of characteristic zero is nilpotent if and only if it admits an invertible left Leibniz-derivation
\end{Th}
\textbf{Proof.} Necessity follows from Lemma \ref{invert_der_exist}, and sufficiency follows from Propositions \ref{Malcev_solv} and \ref{Malcev_nilp}. The Theorem is now proved.

\section{Moens' Theorem for some nonassociative algebras}

\subsection{Jordan algebras}
A commutative algebra $J$ is called \textit{Jordan} if it satisfies the \textit{Jordan identity}:
\begin{eqnarray*} (x^2,y,x) = 0. \end{eqnarray*}
Jordan algebras appeared as one interesting tool for studies in quantum mechanic
in the paper of Jordan, von Neumann and Wigner \cite{jnw34} and
have well known relations with associative, alternative and Lie algebras.
They were  discussed, for example, in \cite{kolca,McCrimmon} and other works.

The classical example of Jordan algebra is constructed as follows: having an associative algebra $A,$ we introduce a new multiplication $\circ$ on its underlying vector space by $x\circ y = \frac{1}{2}(xy + yx)$, and denote the resulting algebra by $A^{(+)}.$ If for Jordan algebra $J$ there exists an associative algebra $A$ such that $J \subseteq A^{(+)},$ then $J$ is called \textit{special}.

If $J$ is a  Jordan algebra over the field of characteristic 0 then $J$ is a direct sum of its solvable radical $Rad(J)$ and a semisimple algebra $S$ isomorphic to $J/Rad(J)$ (Wedderburn decomposition, Albert--Penico--Taft theorem). It is also widely known that $Rad(J)$ coincides with the nilpotent radical of $J$.

Using the well-known fact that a semisimple   Jordan algebra is the direct sum of simple Jordan algebras and the classification of simple  Jordan algebras we can establish that $S$ has a unit.

Since the conditions of Theorem \ref{rad_inv} are true for any Jordan algebra, we have the following theorem:

\begin{Th}
\label{rad_inv_Jord}
Let $J$ be a  Jordan algebra over the field of characteristic 0, $Rad(J)$ be its radical, and $d$ be a left Leibniz-derivation of $J$. Then $d(Rad(J)) \subseteq Rad(J).$
\end{Th}

\medskip

Now we are ready to prove the analogue of Moens' theorem for the class of Jordan algebras:

\begin{Th}
\label{Moens_th_Jord}
A Jordan algebra over a field of characteristic zero is nilpotent if and only if it admits an invertible left Leibniz--derivation.
\end{Th}
\textbf{Proof.}
Necessity follows from Lemma \ref{invert_der_exist} . Now, let $J$ be a Jordan algebra and $d$ be an invertible left Leibniz-derivation of $J.$ Suppose that the semisimple part $S$ of $J$ is nonzero. From Theorem \ref{rad_inv_Jord} it follows immediately that $d(S)=S,$ which contradicts the unitality of $S,$ since for the unit element $1_S$ of $S$ we must have $d(1_S) = nd(1_S),$ hence, in this case $d$ must have nontrivial kernel. The theorem is now proved.

\subsection{$(-1,1)$-algebras.}
An algebra $A$ is called a $(-1,1)$-algebra if it satisfies the following identities:
$$(xy)y = xy^2, (x,y,z) + (z,x,y) + (y,z,x) = 0,$$
where $(x,y,z) = (xy)z - x(yz)$ is the associator of elements $x, y, z \in A.$

The first of the identities above defines the variety of right alternative algebras. In other words, a right alternative algebra $A$ is called $(-1,1)$-algebra, if its adjoint algebra $A^{(-)}$ is a Lie algebra. $(-1,1)$-algebras were studied in many works, such as \cite{Hentz,Mikh} and many others.
 In the paper \cite{Hentz} Hentzel showed that for any  $(-1,1)$-algebra $A$ its right nilpotent radical and nilpotent radical coincide. We will denote them by $Rad(A).$ In \cite{Mikh} it was proved that the Wedderburn principal theorem holds for the class of  $(-1,1)$-algebras, that is, it was proved that any  $(-1,1)$-algebra $A$ is a direct sum of its radical $Rad(A)$ and semisimple part $S$, which is an associative algebra with unit. That implies that arbitrary $(-1,1)$-algebra satisfies the conditions $1)$ and $2)$ of the Theorem \ref{rad_inv}. Now, let $A$ be a  associative algebra and $R$ be its right nilpotent right ideal. Then $R$ is obviously a right nil-ideal of $A$, which implies that $R \subseteq J(A)$, the Jacobson radical of $A$, which coincides with the right nilpotent radical of $A$ in the  case. This and the fact that a semisimple $(-1,1)$-algebra is associative imply that the condition $3)$ of the Theorem 12 holds for arbitrary  $(-1,1)$-algebra. That means that we have proved the following theorem:

\begin{Th}
\label{rad_inv_11}
Let $A$ be a  $(-1,1)$-algebra over the field of characteristic 0, $Rad(A)$ be its radical, and $d$ be a left Leibniz-derivation of $A$. Then $d(Rad(A)) \subseteq Rad(A).$
\end{Th}

\medskip

The proof of this theorem is same as of Theorem \ref{Moens_th_Jord}:
\medskip

\begin{Th}
\label{Moens_th_11}
A $(-1,1)$-algebra over a field of characteristic zero is nilpotent if and only if it admits an invertible left Leibniz--derivation.
\end{Th}

\subsection{Right alternative algebras.}
An algebra $A$ is called \textit{right alternative} if $A$ satisfies the following identity:
\begin{equation*} (x,y,y) = 0. 
\end{equation*}
Right alternative algebras were introduced in \cite{Albert49}.
The class of right alternative algebras is a wide extension of alternative algebras.

It is interesting, that unlike algebras of many well-studied classes, such as Jordan, alternative, Lie and so on, right nilpotent right alternative algebra need not be nilpotent.
The corresponding example of a five-dimensional right nilpotent but not nilpotent algebra belongs to Dorofeev \cite{Dorofeev}.
Its basis is $\{a, b, c , d , e \}$ , and the multiplication is given by the table (zero products of basis vectors are omitted):
$$ab = -ba = ae = -ea = db = -bd = -c, ac = d, bc = e.$$
It is easy to see that the algebra $A$ is not nilpotent. However, it is right nilpotent of index 3, therefore, the identity mapping is an invertible left Leibniz-derivation of order 3.\\
Moreover, one can compute the derivation algebra of this algebra. For any basis element $x \in \{a,b,c,d,e\} = B $ and any linear mapping $\varphi$ we denote $\varphi(x) = \sum_{y \in B} \alpha_{xy}y$. Then any derivation of $A$ has the following matrix form:
$$\varphi =
\begin{pmatrix}
\alpha_{aa} & \alpha_{ba} & 0 & 0 & 0\\
\alpha_{ab} & -\alpha_{aa} & 0 & 0 & 0\\
0 & 0 & \alpha_{ad} + \alpha_{be} & 0 & 0\\
\alpha_{ad} & \alpha_{bd} & 0 & \alpha_{aa} + \alpha_{ad} + \alpha_{be} & \alpha_{ba}\\
\alpha_{ae} & \alpha_{be} & 0 & \alpha_{ab} & -\alpha_{aa} +\alpha_{ad} + \alpha_{bc}
\end{pmatrix},$$
where  $\alpha_{aa}, \alpha_{ab}, \alpha_{ad}, \alpha_{ae}, \alpha_{ba}, \alpha_{bd}, \alpha_{be} \in F.$

It is now easy to see that $A$ possesses invertible derivations. For example, a linear mapping $\varphi$ given by
$$\varphi(a) = a + d, \varphi(b) = -b + e, \varphi(c) = 2c, \varphi(d) = 3d, \varphi(e) = e$$
is an invertible derivation of $A$. Hence, we cannot hope to prove the full analogue of Moens' theorem for right alternative algebras.

\medskip

However, many theorems about right alternative algebras deal more with right nilpotency rahter than nilpotency, so it appears that for right alternative algebras the concept of right nilpotency is more natural.\\
We give a sufficient condition for the right nilpotency of right alternative algebra:

\begin{Th}
\label{Moens_th_right_alt}
A right alternative algebra over a field of characteristic zero  admitting an invertible Leibniz-derivation is right nilpotent.
\end{Th}
\textbf{Proof.}
Let $A$ be a right alternative algebra over field $F$ and $d$ be an invertible Leibniz-derivation of $A.$ Consider the adjoint Jordan algebra $A^{(+)}.$ Expanding Jordan products, applying $d$ to usual $n$-ary products and collecting the terms again, one can show that $d$ is a Leibniz-derivation of $A^{(+)}$ of the same order. Therefore, by Theorem \ref{Moens_th_Jord}, $A^{(+)}$ is a nilpotent Jordan algebra.  
Now, by \cite{skosyr79} we have that $A$ is right nilpotent.
The Theorem is now proved.

\medskip

\subsection{Noncommutative Jordan algebras.}
An algebra $U$ is called a \emph{noncommutative Jordan algebra} if it satisfies the following identities:
$$(x^2y)x=x^2(yx), (xy)x=x(yx).$$
One can see that the class of noncommutative Jordan algebras coincides with the class of algebras in which for arbitrary element $x$ operators $L_x, R_x, L_{x^2}, R_{x^2}$ commute pairwisely.
The class of noncommutative Jordan algebras is rather vast: it includes
all associative algebras, alternative algebras, Jordan algebras, quadratic flexible algebras and anticommutative algebras.

\medskip

Noncommutative Jordan algebras and superalgebras were studied in \cite{McC} and many other works.
It is well known that $U$ is a noncommutative Jordan algebra if and only if $U$ is a flexible algebra such that its symmetrized algebra $U^{(+)}$ is a Jordan algebra.
An algebra $U$ is called \emph{Malcev-admissible}, if its adjoint algebra $U^{(-)}$ is a Malcev algebra. We can also characterize nilpotent noncommutative Jordan Malcev-admissible algebras by Leibniz-derivations:

\begin{Th}
A noncommutative Jordan Malcev-admissible algebra over a field of characteristic zero is nilpotent if and only if it has an invertible Leibniz-derivation.
\end{Th}
\textbf{Proof.} The necessity follows from Proposition \ref{nilp_der}. Now, let $U$ be a noncommutative Jordan algebra such that $U^{(-)}$ is a Malcev algebra, and $d$ be an invertible left Leibniz-derivation of order $n$ of $U.$ It is easy to see that $d$ is an invertible Leibniz-derivation of the same order of both $U^{(+)}$ and $U^{(-)}$. From theorems \ref{Moens_th_Jord} and \ref{Moens_th_Malc} it follows that $U^{(+)}$ and $U^{(-)}$ are nilpotent. From \cite{Skos} it follows that $M(U) = M(U^{(+)})M(U^{(-)}) = M(U^{(-)})M(U^{(+)}),$ hence $M(U)$ is nilpotent and $U$ is nilpotent. The Theorem is now proved.

\subsection{Quasiassociative and quasialternative algebras.}
Some noncommutative Jordan Malcev-admissible algebras can be easily constructed from associative and alternative algebras using the process of mutation.
\medskip

Let $U = (U, \cdot)$ be an algebra over the field $F$ and $\lambda \in F.$ By $U^{(\lambda)}$ we denote an algebra $(U, \cdot_\lambda),$ where
$$x\cdot_\lambda y = \lambda x \cdot y + (1-\lambda)y \cdot x.$$
The algebra $U^{(\lambda)}$ is called a \emph{$\lambda$-mutation} of $U$.
A mutation $A^{(\lambda)}$ of an associative (alternative) algebra $A$ is called a \emph{quasiassociative algebra} (\emph{quasialternative algebra}). 
Quasiassociative and quasialternative algebras were considered in \cite{McC, Skos, Ded} and many other works.

\medskip

We now give the criterion for the nilpotency of quasialternative (in particularly, quasiassociative) algebras:

\begin{Th}
A  quasialternative algebra over a field of characteristic zero is nilpotent if and only if has an invertible Leibniz-derivation.
\end{Th}

\medskip

\subsection{Zinbiel algebras.}
An algebra $A$ is called a Zinbiel  algebra if it satisfies
the following identity:
$$(ab)c=a(bc+cb).$$
As was proved in \cite{dzhuma05},
every finite-dimensional Zinbiel algebra over a field with characteristic zero is nilpotent.
It is easy to see that, the Moens´ theorem is true in the variety of Zinbiel algebras:

\begin{Th}
A Zinbiel algebra over a field of characteristic zero is nilpotent if and only if has an invertible left Leibniz-derivation.
\end{Th}

\section{Appendix: Remarks and open problems}

\subsection{Filippov ($n$-Lie)  algebras.}
An anticommutative $n$-ary algebra $A$ is called \textit{Filippov} ($n$-Lie) if it satisfies the \textit{$n$-ary Jacoby identity}:
\begin{eqnarray*}
[[x_1, \ldots, x_n], y_2, \ldots, y_n]  = \sum [x_1, \ldots, [x_i, y_2, \ldots, y_n], \ldots, x_n]. \end{eqnarray*}
Filippov algebras appeared in the paper of Filippov \cite{FIL85} and were used as an interesting tool for the study of mathematical physics
(see, for example,  survey \cite{AZI}).

Williams proved that the Theorem of Jacobson is not true in the case of Filippov algebras.
Namely, he constructed a non-simple non-nilpotent Filippov algebra with an invertible derivation \cite[Example 3.1]{Williams}.
For $n$-dimensional vector space with the basis $x_1, \ldots, x_n$ we define the anticommutative multiplication $[x_1, \ldots, x_n]=x_2.$
He made a mistake in the case $n$ is even.
For example, if $n=4$ Williams' mapping is not a derivation.
We can correct their example:
$$D(x_1)=x_1, D(x_2)=x_2, D(x_3)=-\frac{1}{2}x_3, D(x_4)=-\frac{1}{2}x_4, D(x_j)=(-1)^jx_j, \mbox{ for }j>4.$$

The main example of Filippov algebras is a simple $n$-Lie algebra $D_{n+1}$ in dimension $n+1$.
On the basis elements $e_1, \dots, e_{n+1}$ of $D_{n+1}$ we define the product the following way:
$$ [e_1, \dots, \hat{e_i}, \dots, e_{n+1}] = (-1)^{n+i+1}e_i.$$
As follows from \cite{KayFil},
for any $n>2$, the simple Filippov algebra $D_{n+1}$ over a field of characteristic zero has an invertible derivation.
For example,
the linear mapping with the matrix
$$\sum( e_{ii+1}-e_{i+1i}), \mbox{  if }n+1=2k,$$
or
$$\sum( e_{ii+1}-e_{i+1i}) +e_{1n+1}-e_{n+11}, \mbox{  if }n=2k,$$
is an invertible derivation of $D_{n+1}.$
\subsection{$n$-Ary Malcev algebras.}
An anticommutative $n$-ary algebra $A$ is called \textit{$n$-ary Malcev} algebra
 if it satisfies the \textit{$n$-ary Malcev identity}:
\begin{eqnarray*}
−J(zR_x,x_2, \ldots ,x_n;y_2, \ldots ,y_n) = J(z,x_2, \ldots ,x_n;y_2, \ldots ,y_n)R_x,\end{eqnarray*}
where
$$J(x_1,x_2, \ldots ,x_n;y_2,...,y_n)=
[[x_1, \ldots, x_n], y_2, \ldots, y_n]  - \sum [x_1, \ldots, [x_i, y_2, \ldots, y_n], \ldots, x_n],
R_x=R_{x_2, \ldots, x_n}.$$

$n$-Ary Malcev algebras appeared in the paper of Pojidaev \cite{POJ01} as one generalization of Filippov and Malcev algebras.
As follows from \cite{POJSAR},
the simple
$8$-dimensional  ternary Malcev algebra over a field of characteristic zero admits an invertible derivation.
For example,
the linear mapping with the matrix $$e_{23}-e_{32}-e_{14}+e_{41}+e_{56}-e_{65}+e_{78}-e_{87}.$$

\subsection{Remarks.}
We finish this paper with some remarks.

Note that the theorem only holds for characteristic zero:

\begin{Remark}
For an arbitrary algebra over a field of positive characteristic $p$, the identity map is a Leibniz-derivation of order $p+1$.
\end{Remark}

We also note that the theorem fails for infinite-dimensional algebras:

\begin{Remark}
The free algebra of any homogeneous variety $\Omega$ not consisting entirely of nilpotent algebras admits an invertible derivation, but it is not nilpotent.
\end{Remark}

\medskip Indeed, it suffices to consider a derivation which acts identically on the generators of the algebra.

\subsection{Open problems.}
We now know that an analogue of the theorem of Moens is true for
Lie, Leibniz, Zinbiel, associative, alternative, Jordan, Malcev, $(1,-1)$-, right alternative,
quasiassociative, quasialternative and noncommutative Jordan Malcev-admissible algebras;
but not true for Filippov and $n$-ary Malcev algebras.
Therefore, we can state some problems which, as far as the authors know, are still open:

\

{\bf Problem 1.} Does the Moens' theorem hold in the classes of:

(a) Binary Lie,  Novikov, right-symmetric, bi-commutative and other algebras;

(b) Alternative, Jordan, Lie, Malcev and other  superalgebras?

\

{\bf Problem 2.} Is there a variety of $n$-ary algebras $(n>2)$,
where an analogue of Theorem of Jacobson (and, more generally, Moens') holds?

\medskip



\end{document}